\documentclass[10pt, a4paper]{article}
\usepackage{latexsym}
\usepackage{indentfirst}
\usepackage{amsmath}
\usepackage{amssymb}

\begin{document}

\title{Errata to "An example of Bautin-type bifurcation in a delay differential
equation", JMAA, 329(2007), 777-789}
\author{Anca-Veronica Ion\\
"Gh. Mihoc-C. Iacob" Institute of Mathematical Statistics\\
and Applied Mathematics of the Romanian Academy, \\
13, Calea 13 Septembrie, Bucharest, Romania \\
anca-veronica.ion@ima.ro}

\date{}
\maketitle

\begin{abstract}
Some errors contained in the author's previous article "An example of
Bautin-type bifurcation in a delay differential equation", JMAA, 329(2007),
777-789, are listed and corrected.
\end{abstract}

In our work \cite{AVI}, we considered the delay differential equation,
\begin{equation}\label{ex}\dot{x}=ax(t-r)+x^2(t)+cx(t)x(t-r),
\end{equation}
and looked for values of the parameters where the conditions for the occurrence of Bautin type bifurcation around the equilibrium point $x=0$ are fulfilled.

For this  we first proved that for the linearized around
$x(t)=0$ equation, at $$a=-1,\,\, r=\pi/2,$$  two eigenvalues
$\lambda_{1,2}=\pm i$ exist,  while all the other eigenvalues have
negative real part. Thus we were entitled to consider the reduction
of the problem to the two-dimensional center manifold for these
values of the parameters. The reduced problem is a two-dimensional
system of differential equations, that can be written as an ODE for a complex valued function
\begin{equation}\label{red-eq}\dot{z}=\pm i z + \sum_{j+k\geq 2}\frac{1}{j!k!}g_{jk}(c)z^{j}\overline{z}^{k}.
\end{equation}

 For such problems the Bautin
bifurcation was studied in \cite{K} and we followed the method therein for our study. In order to find Bautin
bifurcation points we computed the first Lyapunov coefficient and
found that this is zero for
\[c_{1,2}=\frac{18-7\pi\pm\sqrt{36+212\pi+\pi^2}}{2(3\pi-2)}.
\]

In this note we intend to correct two distinct type of mistakes that, we, unhappilly, made in \cite{AVI}.

\vspace{0.3cm}

\textbf{I.}  To determine whether the bifurcation point presents a higher order
degeneracy or is a proper Bautin bifurcation point, we computed
the second Lyapunov coefficient for the reduced on the center manifold problem. For this we needed
$w_{21}(0),\,w_{21}(-r),$ where $w_{21}(\cdot)\in \mathcal{C}([-r,0],\,\mathbb{R})$ is a coefficient
of the series of powers of the function whose graph is the center
manifold ($w(z,\overline{z})(\cdot)=\sum_{j+k\geq 2}\frac{1}{j!k!}w_{jk}(\cdot)z^{j}\overline{z}^{k}$ see \cite{AVI}).

The two algebraic equations that yield $w_{21}(0)$ and $w_{21}(-r)$ proved to be  dependent,
and at that moment we have chosen arbitrarily $w_{21}(0)=0$ and we
computed $w_{21}(-r)$ from one of the two equations. This is a
mistake and we want to correct it here.

By studying more carefully the problem of computing $w_{21}(0)$ and
$w_{21}(-r)$, we found out that these can be uniquely determined by
using a perturbation technique. This result was published in
\cite{AVI1}.

The formula obtained there for $w_{21}(0),$ adapted to problem \eqref{ex}, is
\begin{equation}\label{w21}w_{21}(0)=
\frac{f_{21}\langle \Psi_{1}+\Psi_{2}, \,\rho
\rangle-2g_{11}\langle \widetilde{\rho}, w_{20}\rangle
-(g_{20}+2\overline{g}_{11})\langle \widetilde{\rho},
w_{11}\rangle -\overline{g}_{02}\langle \widetilde{\rho},
w_{02}\rangle }{2r i +2}.
\end{equation}
Here $\Psi_{1}(\zeta)=2\frac{2-\pi i}{4+\pi^2}e^{- i\zeta}, \;\zeta\in [0,r],$
$\Psi_{2}=\overline{\Psi}_{1}$, $f_{21}=g_{21}/\Psi_{1}(0),$  $g_{ij}$ are the coefficients of \eqref{red-eq}$, \rho(s)=-2se^{ is}, \;s\in
[-r,0]$ and $\widetilde{\rho}(\zeta)=-2\zeta e^{-
i\zeta}, \;\zeta\in [0,r]$,  while by the brackets
$\langle\,\cdot\,,\,\cdot\,\rangle$ we denote the bilinear form
defined in the study of delay differential equations (see
\cite{AVI1} and the references therein).

By using formula \eqref{w21} we
found,  in the case of $c_1(\approx1.52799)$:
\[w_{21}(0)=0.4748-0.4547
i,\,w_{21}(-r)=1.4926-1.9467 i,\,l_2=1.305.\]

Hence, by the theory concerning the Bautin bifurcation, in the
parameters plane, in a neighborhood of the point $a=-1,\,c=c_1,$
there is a zone where an unstable manifold exists and for parameters
$a,c$ in a subset of this zone, two periodic orbits (one inside the
other) exist on the unstable manifold. The inner periodic (closed)
orbit is attracting, while the outer one is repelling.

We then analyzed the case of $c_2(\approx-2.06554)$ and found:
\[w_{21}(0)=-0.2687-0.0084i,\,w_{21}(-r)=-4.1734-1.7929i,\,l_2=10.421.\]

This shows that equation \eqref{ex} presents the same type of Bautin
type bifurcation for both pairs of parameters $a=-1,\,c=c_{1}$ and $a=-1,\, c=c_2$.

\vspace{0.5cm}

\textbf{II.} We also noticed some other errors in \cite{AVI}, that
we correct  here:

\begin{enumerate}
  \item at pg. 8 (784 in JMAA), $w_{20}(0,c)$ should be $$w_{20}(0,c)=F_{20}\left[\frac{4(\pi+ 4i)}{3(4+\pi^2)}-\frac{1+2i}{5}\right]=2(1-ic)\left[\frac{4(\pi+ 4i)}{3(4+\pi^2)}-\frac{1+2i}{5}\right];$$
  \item at pg. 10 (787 in JMAA), $F_{31}$  should be  \[F_{31}=c_1[3w_{21}(-r)+w_{30}(-r)+iw_{30}(0)-3iw_{21}(0)+\]
      \[+3w_{20}(0)w_{11}(-r)+3w_{11}(0)w_{20}(-r)]+6w_{11}(0)w_{20}(0)+6w_{21}(0)+2w_{30}(0);\]
  \item  at pg. 10 (787 in JMAA), $F_{22}$ should be
  \[F_{22}=c_1[2w_{12}(-r)+2w_{21}(-r)+w_{20}(0)w_{02}(-r)+4w_{11}(0)w_{11}(-r)+w_{02}(0)w_{20}(-r)
  +\]
  \[+2iw_{21}(0)-2iw_{12}(0)]+2w_{20}(0)w_{02}(0)+4w_{11}(0)^2+4w_{12}(0)+4w_{21}(0).
  \]
\end{enumerate}

We apologize to the readers of \textit{Journal of Mathematical Analysis and Applications } for the errors listed and corrected above.

\newpage

\setcounter{page}{1}

\begin{center}
\Large{An example of Bautin-type bifurcation in a delay differential equation}\footnote{published in Journal of Mathematical Analysis and Appplication, 329(2007), 777-789}\\

\vspace{0.5cm}

Anca-Veronica Ion\\ 

\vspace{0.3cm}

\small{University of Pite\c sti, Faculty of Mathematics and Computer Sciences,\\
1, str. T\^argu din Vale, Pite\c sti, 110040, Arge\c s,  Romania;\\
  e-mail: averionro@yahoo.com}
\end{center}

\begin{abstract}
In a previous paper we gave sufficient conditions for a system of delay
differential equations to present Bautin-type bifurcation. In the present
work we present an example of delay equation that satisfies these conditions.
\end{abstract}

\section{\protect\bigskip Introduction}

\setcounter{equation}{0}

\bigskip In \cite{I} the system of delay differential equations
\begin{eqnarray}
\overset{\cdot }{x}(t) &=&A\left( \alpha \right) x\left( t\right) +B\left(
\alpha \right) x\left( t-r\right) +f\left( x\left( t\right) ,x\left(
t-r\right) ,\alpha \right) ,  \label{pr} \\
x\left( s\right) &=&\phi \left( s\right) ,\;\;\;s\in \left[ -r,0\right] ,
\notag
\end{eqnarray}
with $x\left( t\right) =\left( x_{1}\left( t\right) ,...,x_{n}\left(
t\right) \right) \in \mathbb{R}^{n},\;\alpha =\left( \alpha _{1},\alpha
_{2}\right) \in \mathbb{R}^{2},\;A\left( \alpha \right) ,\;B\left( \alpha
\right) \;n\times n$ real matrices is considered. Here $f=\left(
f_{1},...f_{n}\right) $ is continuously differentiable on its domain of
existence, $D\subset \mathbb{R}^{2n+2}.$ It is also assumed that $f\left(
0,0,\alpha \right) =0$ and the differential of $f$ in the first two
vectorial variables, calculated at $\left( 0,0,\alpha \right) $ is equal to
zero. $\phi $ belongs to the Banach space $C\left( \left[ -r,0\right] ,%
\mathbb{R}^{n}\right) .$

For this system we give in \cite{I} a theorem providing sufficient
conditions for the appearance of Bautin-type bifurcation.

Bautin bifurcations are degenerated Hopf bifurcations. As it is known, \cite
{K}, for two-dimensional systems of ODEs depending on a scalar parameter $%
\alpha $, Hopf bifurcation around a branch of equilibrium points appears
when there is a certain value of the parameter, $\alpha _{0}$, at which:

- a pair of purely imaginary eigenvalues of the linear part exists,

- the real part of the eigenvalues (that is zero at $\alpha _{0}$) has
non-zero derivative at $\alpha _{0}$,

- the first Lyapunov coefficient at $\alpha _{0}\;$is non-zero.

The first Lyapunov coefficient is a number defined as follows. The
two-dimensional system of real equations is written as a single complex
equation
\begin{equation*}
\overset{\cdot }{z}=\lambda z+g(z,\overset{\_}{z},\alpha ),
\end{equation*}
and the first Lyapunov coefficient, $l_{1}(\alpha )$ is defined in terms of
the coefficients (up to the third degree) of the series

\begin{equation*}
g\left( z,\overset{\_}{z},\alpha \right) =\sum_{j+k\geq 2}\frac{1}{j!k!}%
g_{jk}\left( \alpha \right) z^{j}\overset{\_}{z}^{k},
\end{equation*}
(see (\ref{l1}) below).

Now, also for an ODEs system, let us assume that the parameter $\alpha $ is
bidimensional. When $l_{1}(\alpha _{0})=0$ and a second Lyapunov
coefficient, $l_{2}(\alpha )$ (that is defined in terms of the coefficients
up to the fifth degree terms of the above series - see Section 6.2) is
non-zero at $\alpha _{0},$ the Bautin bifurcation takes place \cite{K}. It
is characterized by the appearance, for the parameters in a neighborhood of $%
\alpha _{0}$, of two limit cycles (one inside the other).

By using the reduction of the problem (\ref{pr}) to its central manifold, we
extended in the main theorem in \cite{I} the above ideas to the class of
systems of delay differential equations (\ref{pr}).

The problem that arised after this theorem was proved, was whether there is
any delay equation that satisfies its hypotheses, or not.

In this paper we present a scalar differential delay equation that satisfies
the hypotheses of our theorem.

We set below the theoretical frame we used and our result, as a starting
point for the rest of the paper.

\section{Theoretical framework}

\bigskip Let us consider the solutions of the equation
\begin{equation*}
\det \left( \lambda I-A\left( \alpha \right) -e^{\lambda r}B\left( \alpha
\right) \right) =0,
\end{equation*}
where $I$ is the $n$-dimensional unity matrix. These are the eigenvalues of
the infinitesimal generator of the linearized problem obtained from (\ref{pr}%
) (see \cite{H}, \cite{HL}). Let us consider the hypothesis:

\textit{H1. There is an open set U in the parameter plane such that for
every }$\alpha \in U$\textit{, there is a pair of complex conjugated simple
eigenvalues }$\lambda _{1,2}\left( \alpha \right) =\mu \left( \alpha \right)
$\textit{\ }$\pm i\omega \left( \alpha \right) $, \textit{with the property
that there is a }$\alpha _{0}\in U\;\,$\textit{such that }$\lambda
_{1,2}\left( \alpha _{0}\right) =\pm i\omega \left( \alpha _{0}\right) =\pm
i\omega _{0}$\textit{, with }$\omega _{0}>0$\textit{\ and there is an }$%
\varepsilon >0\;$\textit{\thinspace such that for every }$\alpha \in U,$ $%
\mu \left( \alpha \right) >-\varepsilon ,$ \textit{while all other
eigenvalues }$\lambda $\textit{\ have }$\mathrm{Re}\lambda <-\varepsilon $%
\textit{. }

It is important to assume that, as $\alpha $ varies in $U$, $\mu \left(
\alpha \right) $ takes both positive and negative values.This is usually
expressed by the hypothesis $\frac{d\mu }{d\alpha }\left( \alpha _{0}\right)
\neq 0,$ but in our case it will be covered by hypothesis \textit{H2 }below.

Let $\varphi _{1}\left( \alpha \right) ,\;\varphi _{2}\left( \alpha \right)
(=\overset{\_}{\varphi }_{1}\left( \alpha \right) )\in C\left( \left[ -r,0%
\right] ,\mathbb{R}^{n}\right) $ be the two eigenvectors corresponding to $%
\lambda _{1}\left( \alpha \right) $ respectively $\lambda _{2}\left( \alpha
\right) \;$ (these are simple eigenvalues). Let also $\mathbb{M}_{\left\{
\lambda _{1,2}\left( \alpha \right) \right\} }$ be the space spanned by $%
\varphi _{1}\left( \alpha \right) ,\;\varphi _{2}\left( \alpha \right) $.

For the values of $\alpha \in U$ such that $\mu \left( \alpha \right) >0$
there is a two-dimensional local invariant manifold, the unstable manifold
of the equilibrium point $0$ (see \cite{H}, \cite{HL}, \cite{MNO}). \ For $%
\alpha _{0}$ there is a two-dimensional local central manifold. In both
cases the manifold is the graph of a differentiable application $w_{\alpha }$
defined on $\mathbb{M}_{\left\{ \lambda _{1,2}\left( \alpha \right) \right\}
}.$ We denote the local invariant manifold by $W_{loc}\left( \alpha \right)
. $

The restriction of the equation (\ref{pr}) to the invariant manifold for the
values of $\alpha $ mentioned above is
\begin{equation}
\overset{\cdot }{z}\left( t\right) =\lambda _{1}\left( \alpha \right)
z\left( t\right) +\psi _{1}\left( \alpha \right) \left( 0\right) f\left( %
\left[ S_{\alpha }\left( t\right) \phi \right] \left( 0\right) ,\left[
S_{\alpha }\left( t\right) \phi \right] \left( -r\right) ,\alpha \right) ,
\label{ecz}
\end{equation}
where $\psi _{1}\left( \alpha \right) $ is a certain eigenvector of the
adjoint problem (\cite{H}, \cite{HMO}).

If we take $\phi \in W_{loc}\left( \alpha \right) ,$ then $S_{\alpha }\left(
t\right) \phi \in W_{loc}\left( \alpha \right) $ and thus

\begin{center}
\begin{equation}
S_{\alpha }\left( t\right) \phi \left( s\right) =z\left( t\right) \varphi
_{1}\left( s\right) +\overset{\_}{z}\left( t\right) \overset{\_}{\varphi }%
_{1}\left( s\right) +w_{\alpha }\left( s,z\left( t\right) ,\overset{\_}{z}%
\left( t\right) \right) .  \label{Salfa}
\end{equation}
\end{center}

This implies that $f\left( \left[ S_{\alpha }\left( t\right) \phi \right]
\left( s\right) ,\left[ S_{\alpha }\left( t\right) \phi \right] \left(
s-r\right) ,\alpha \right) $ is a function of $z,\overset{\_}{z}$ and it can
be written as a series of powers as

\begin{center}
\begin{equation}
f\left( S_{\alpha }\left( t\right) \phi \left( s\right) ,S_{\alpha }\left(
t\right) \phi \left( s-r\right) ,\alpha \right) =\sum\limits_{j+k\geq 2}%
\frac{1}{j!k!}F_{jk}\left( s,\alpha \right) z^{j}\overset{\_}{z}^{k}.
\label{F-z}
\end{equation}
\end{center}

$\bigskip $Then we can write $\psi _{1}\left( \alpha \right) \left( 0\right)
f\left( \left[ S_{\alpha }\left( t\right) \phi \right] \left( 0\right) ,%
\left[ S_{\alpha }\left( t\right) \phi \right] \left( -r\right) ,\alpha
\right) \;$also as a function of $z\left( t\right) ,\overset{\_}{z}\left(
t\right) ,$ namely
\begin{equation}
\psi _{1}\left( \alpha \right) \left( 0\right) f\left( \left[ S_{\alpha
}\left( t\right) \phi \right] \left( 0\right) ,\left[ S_{\alpha }\left(
t\right) \phi \right] \left( -r\right) ,\alpha \right) =g\left( z\left(
t\right) ,\overset{\_}{z}\left( t\right) ,\alpha \right) .  \label{rel f-g}
\end{equation}

and:
\begin{equation*}
g\left( z\left( t\right) ,\overset{\_}{z}\left( t\right) ,\alpha \right)
=\sum_{j+k\geq 2}\frac{1}{j!k!}g_{jk}\left( \alpha \right) z\left( t\right)
^{j}\overset{\_}{z}\left( t\right) ^{k}.
\end{equation*}

Equation (\ref{ecz}) becomes
\begin{equation*}
\overset{\cdot }{z}\left( t\right) =\lambda _{1}\left( \alpha \right)
z\left( t\right) +\sum_{j+k\geq 2}\frac{1}{j!k!}g_{jk}\left( \alpha \right)
z\left( t\right) ^{j}\overset{\_}{z}\left( t\right) ^{k}.
\end{equation*}

For this equation we can study Bautin bifurcation as in \cite{K}. We
consider the first and second Lyapunov coefficients defined in \cite{K},
that are functions of $g_{ij}$. We remind that
\begin{equation}
l_{1}\left( \alpha _{0}\right) =\frac{1}{2\omega _{0}^{2}}\mathrm{Re}\left(
ig_{20}\left( \alpha _{0}\right) g_{11}\left( \alpha _{0}\right) +\omega
_{0}g_{21}\left( \alpha _{0}\right) \right) ,  \label{l1}
\end{equation}
while $l_{2}\left( \alpha _{0}\right) $ is a much more complicated
expression.

We also define
\begin{equation*}
\nu _{1}=\frac{\mu \left( \alpha \right) }{\omega \left( \alpha \right) }%
,\;\nu _{2}=l_{1}\left( \alpha \right) ,
\end{equation*}
and $\nu =\left( \nu _{1},\nu _{2}\right) .$ Let us consider the following
hypothesis:

\textit{H2. }$l_{1}\left( \alpha _{0}\right) =0,$ $l_{2}\left( \alpha
_{0}\right) >0,$ \textit{and the map }$\left( \alpha _{1},\alpha _{2}\right)
\rightarrow \left( \nu _{1},\nu _{2}\right) \,\ $\textit{is regular at }$%
\alpha _{0}.$

Now we can state \ the main result of \cite{I}.

\textbf{Theorem }\textit{If H1, H2 are satisfied for eq. (\ref{pr}), then at
}$\alpha _{0}\;$\textit{a Bautin-type bifurcation takes place. }

\textit{That is there is a neighbourhood }$U_{1}\;$\textit{of }$\alpha _{0}\;
$\textit{in the }$\alpha \;$\textit{plane having a subset }$V^{\ast }$%
\textit{\ (with }$\alpha _{0}\;$\textit{as a limit point) with the property
that for every }$\alpha \in V^{\ast }$\textit{, the restriction of problem (%
\ref{pr}) to the unstable manifold has two limit cycles (one interior to the
other). }

\section{The scalar equation}

Let us consider the equation
\begin{equation}
x^{\prime }=ax(t-r)+x^{2}\left( t\right) +cx(t)x(t-r),  \label{mainec}
\end{equation}
with $r=\frac{\pi }{2}.$ We will study the equation around the equilibrium
solution $x(t)=0.$

Define $\alpha =\left( a,c\right) .$ The linear part of (\ref{mainec}) is
\begin{equation}
x^{\prime }=ax(t-r).  \label{linear_ec}
\end{equation}

Let us consider the function
\begin{equation}
\eta (s)=\left\{
\begin{array}{c}
-a,\;\;\;\;\;\;\;\;\;\;\;\;s=-r \\
0,\;\;\;\;\;\;\;\;\;\;\;s\in \left( -r,0\right] .
\end{array}
\right.  \label{eta}
\end{equation}

We observe that, by defining
\begin{equation*}
L\varphi =\int_{-r}^{0}\varphi \left( s\right) d\eta \left( s\right) ,
\end{equation*}
and $x_{t}(s)=x(t+s)\;$for $s\in \lbrack -r,0],\;\ $equation \ (\ref
{linear_ec}) may be written as

\begin{center}
$x^{\prime }=Lx_{t}.$
\end{center}

\section{The eigenvalues}

The characteristic equation is $\lambda -ae^{-\lambda r}=0,$ \ where $%
\lambda =\mu +i\omega $ \ and it is equivalent to the system of two equations

\begin{center}
$\mu -ae^{-\mu r}\cos \omega r=0,$

$\omega +ae^{-\mu r}\sin \omega r=0.$
\end{center}

These are equivalent with

\begin{eqnarray}
\omega &=&\pm \sqrt{a^{2}e^{-2\mu r}-\mu ^{2}}.  \label{omega} \\
\cos \sqrt{a^{2}e^{-2\mu r}-\mu ^{2}}r &=&\frac{\mu }{a}e^{\mu r}.
\label{ec_cos}
\end{eqnarray}

We see that at $a_{0}=-1,\;r=\frac{\pi }{2}\;$ the pairs $\omega =\pm
1,\;\mu =0\;$\ are solutions of the above equations.

In order to study equation (\ref{ec_cos}) we define $y=\mu r,$ and obtain
the new equation
\begin{equation}
\cos \sqrt{\frac{a^{2}r^{2}}{e^{2y}}-y^{2}}=\frac{y}{ar}e^{y},
\label{ecuatiey}
\end{equation}
that accepts the solution $y=0$ for $a_{0}=-1.$

\bigskip Hence $\lambda _{1,2}\left( a_{0}\right) =\pm i\omega \left(
a_{0}\right) =\pm i\omega _{0}=\pm i,\;\omega _{0}=1.$

\textbf{Proposition 1. }\textit{There is a open neighborhood }$V_{-1}$
\textit{of }$a_{0}=-1$ \textit{such that for every }$a\in V_{-1}$\textit{,
there is a pair of complex conjugated simple eigenvalues }$\lambda
_{1,2}\left( a\right) =\mu \left( a\right) $\textit{\ }$\pm i\omega \left(
a\right) $ \textit{such that} \textit{for every }$a\in V_{-1},\;\mu (a)>-%
\frac{1}{8},$ \textit{\ and all other eigenvalues }$\lambda $ \textit{have }$%
\mathrm{Re}\lambda <-\frac{1}{8}$\textit{. }

\textbf{Proof. }We consider the function$\;G(a,y)=\cos \sqrt{\frac{a^{2}r^{2}%
}{e^{2y}}-y^{2}}-\frac{y}{ar}e^{y},$ and \ we observe that $G\left(
a_{0},0\right) =0$ and

\begin{center}
$\frac{\partial G}{\partial y}=\frac{\frac{a^{2}r^{2}}{e^{2y}}+y}{\sqrt{%
\frac{a^{2}r^{2}}{e^{2y}}-y^{2}}}\sin \sqrt{\frac{a^{2}r^{2}}{e^{2y}}-y^{2}}-%
\frac{y+1}{ar}e^{y},$

$\frac{\partial G}{\partial y}\left( -1,0\right) =\frac{\pi }{2}+\frac{2}{%
\pi }>0.$
\end{center}

The implicit functions theorem implies the existence of: a neighborhood \ $%
W_{-1}$ of $-1$, a neighborhood $W_{0}$ of $0$ and an unique function $%
y:W_{-1}\rightarrow W_{0}$ such that $G\left( a,y(a)\right) =0$ for every $%
a\in W_{-1}.$

Thus $\mu \left( a\right) =\frac{1}{r}y(a)$, $\omega \left( a\right) =\pm
\sqrt{\frac{a^{2}}{e^{2y\left( a\right) }}-\frac{y\left( a\right) ^{2}}{r^{2}%
}}$ and we can define the eigenvalues
\begin{equation*}
\lambda _{1,2}\left( a\right) =\frac{1}{r}y\left( a\right) \pm i\sqrt{\frac{%
a^{2}}{e^{2y\left( a\right) }}-\frac{y\left( a\right) ^{2}}{r^{2}}}.
\end{equation*}

For $a_{0}=-1,\;y\left( a_{0}\right) =0,\;$we have $\sqrt{\frac{a^{2}r^{2}}{%
e^{2y}}-y^{2}}=\frac{\pi }{2}.\;$Let us denote by $m$ a positive integer
such that
\begin{equation*}
\frac{\frac{1}{4}r^{2}}{e^{2\frac{\pi }{m}}}-\left( \frac{\pi }{m}\right)
^{2}>0,\;e^{\frac{\pi }{m}}<\frac{4}{3}.
\end{equation*}

There is a neighbourhood of $-1,$ $W_{-1}^{m}$ $\subset W_{-1}\;$such that
for $a\in W_{-1}^{m}$, $\ y\left( a\right) \in \left( -\frac{\pi }{m},\frac{%
\pi }{m}\right) .\;$We shall take $V_{-1}^{m}=\left( -\frac{3}{2},-\frac{1}{2%
}\right) \cap W_{-1}^{m}.\;$This implies

\begin{center}
$\sqrt{\frac{\frac{1}{4}r^{2}}{e^{2\frac{\pi }{m}}}-\left( \frac{\pi }{m}%
\right) ^{2}}\;\leq \sqrt{\frac{a^{2}r^{2}}{e^{2y\left( a\right) }}-y\left(
a\right) ^{2}}\leq \sqrt{\frac{9}{4}r^{2}e^{\frac{2\pi }{m}}}=\frac{3\pi }{4}%
e^{\frac{\pi }{m}}<\pi ,$
\end{center}

\noindent for $a\in $ $V_{-1}^{m}.$

We consider the equation

\begin{center}
$\frac{a^{2}r^{2}}{e^{2y}}-y^{2}=0$
\end{center}

and denote by $y_{r}\left( a\right) \;$its positive solution. We look for
solutions of (\ref{ec_cos}) only at the left of $y_{r}\left( a\right) $,
since only there the expression $\sqrt{\frac{a^{2}r^{2}}{e^{2y\left(
a\right) }}-y\left( a\right) ^{2}}\;$is real.

Since the function $u\rightarrow \frac{\sin u}{u}$ is decreasing on $\left[
0,\pi \right] ,\;$for $y\in \left[ 0,y_{r}\left( a\right) \right] ,\;$and $%
a<0\;$we have

\begin{center}
$\frac{\partial G}{\partial y}>\frac{\frac{a^{2}r^{2}}{e^{2y}}+y}{\sqrt{%
\frac{a^{2}r^{2}}{e^{2y}}-y^{2}}}\sin \sqrt{\frac{a^{2}r^{2}}{e^{2y}}-y^{2}}-%
\frac{y+1}{ar}e^{y}\geq -\frac{y+1}{ar}e^{y}>0.$
\end{center}

Thus, there are no solutions at the right of $y\left( a\right) $ for $a\in
V_{-1}^{m}.$

Now, for $y\in \left( -\frac{\pi }{m},0\right) $

\begin{center}
$\sqrt{\frac{\pi ^{2}}{16}-\left( \frac{\pi }{m}\right) ^{2}}\leq \sqrt{%
\frac{a^{2}r^{2}}{e^{2y}}-y^{2}}\leq \frac{\pi }{2}e^{\frac{\pi }{m}}.$
\end{center}

We choose $m=16$ and denote the neighborhood $V_{-1}^{16}$ by $V_{-1}.\;$Then

\begin{center}
$\sqrt{\frac{\pi ^{2}}{16}-\left( \frac{\pi }{16}\right) ^{2}}\;=\frac{\sqrt{%
15}\pi }{8}\geq \frac{\pi }{6},$

$\frac{\pi }{6}\leq \sqrt{\frac{a^{2}r^{2}}{e^{2y}}-y^{2}}\leq \frac{2\pi }{3%
}$
\end{center}

and thus

\begin{center}
$\frac{\partial G}{\partial y}\left( a,y\right) =\frac{\frac{a^{2}r^{2}}{%
e^{2y}}+y}{\sqrt{\frac{a^{2}r^{2}}{e^{2y}}-y^{2}}}\sin \sqrt{\frac{a^{2}r^{2}%
}{e^{2y}}-y^{2}}-\frac{y+1}{ar}e^{y}>$

$\;>\left( \frac{\pi ^{2}}{4}+y\right) \frac{\sqrt{3}}{2}\frac{3}{2\pi }-%
\frac{y+1}{ar}e^{y}\geq $

$\;\;\;\;\;\;\;\;\;\geq \left( \frac{\pi ^{2}}{4}-\frac{\pi }{16}\right)
\frac{3\sqrt{3}}{4\pi }+2\frac{-\frac{\pi }{16}+1}{r}e^{y}>0.$
\end{center}

Hence we have no solutions with $\mu \geq -\frac{\pi }{16}r=-\frac{\pi }{16}%
\frac{2}{\pi }=-\frac{1}{8}$ besides $\mu \left( a\right) =\frac{2}{\pi }%
y\left( a\right) .$

To summarize, for each $a\in V_{-1}\;$ we have the following:

\noindent - \ \ there is a pair of eigenvalues of (\ref{ecuatiey}), namely

\begin{center}
$\lambda _{1,2}\left( a\right) =\frac{2}{\pi }y\left( a\right) \pm i\sqrt{%
\frac{a^{2}}{e^{2y\left( a\right) }}-\frac{y\left( a\right) ^{2}}{r^{2}}}$
\end{center}

\noindent with $y(a)$ defined above,

\noindent - $\ \ \mu \left( a\right) >-\frac{1}{8},$

\noindent - \ \ all other eigenvalues $\lambda $\ have $\mathrm{Re}\lambda <-%
\frac{1}{8},$

\noindent - \ \ for $a_{0}=-1,$ $\lambda _{1,2}\left( a_{0}\right) =\pm
i.\square $

It follows that the hypothesis \textit{H1} is satisfied by our equation.

\section{The eigenvectors at $a_{0}$}

The eigenvectors \cite{HL} corresponding to $\lambda _{1,2}\left(
a_{0}\right) $ are $\varphi _{1}\left( s\right) =e^{is},\;\;\varphi
_{2}\left( s\right) =e^{-is},$ $s\in \left[ -r,0\right] ,$ and we denote by $%
\mathbb{M}_{\left\{ \lambda _{1,2}\left( a_{0}\right) \right\} }$ the
eigenspace spanned by them.

The eigenvectors for the adjoint problem are $\phi _{1}\left( s\right)
=e^{-is},\;\phi _{2}\left( s\right) =e^{is},$ $s\in \left[ 0,r\right] $. Let
us denote by $\mathbb{M}_{\left\{ \lambda _{1,2}\left( a_{0}\right) \right\}
}^{\ast }$ the space spanned by $\left\{ \phi _{1},\phi _{2}\right\} $ in $%
C\left( \left[ 0,r\right] ,\mathbb{R}^{n}\right) .$

We define the bilinear form $\left( \chi (.),\varphi (.)\right) :\mathbb{M}%
_{\left\{ \lambda _{1,2}\left( a_{0}\right) \right\} }^{\ast }\times \mathbb{%
M}_{\left\{ \lambda _{1,2}\left( a_{0}\right) \right\} }\rightarrow \mathbb{C%
}$,

$\left( \chi (.),\varphi (.)\right) =\overset{\_}{\chi }(0)\varphi
(0)-\int_{-r}^{0}\int_{0}^{\theta }\chi \left( \xi -\theta \right) \psi
\left( \xi \right) d\xi d\eta \left( \theta \right) $, with $\eta $ defined
by$\;$(\ref{eta}).

Let the numbers $e_{ij}$ be defined by $e_{ij}=\left( \phi _{i}(.),\varphi
_{j}(.)\right) .$ We find that the matrix $E=\left( e_{ij}\right) _{1\leq
i,j\leq 2}$ is

\begin{center}
$E=\left[
\begin{array}{cc}
\frac{2+\pi i}{2} & 0 \\
0 & \frac{2-\pi i}{2}
\end{array}
\right] .$
\end{center}

The vectors $\psi _{1},\psi _{2}$ given by

\begin{center}
\begin{equation}
\left(
\begin{array}{c}
\psi _{1} \\
\psi _{2}
\end{array}
\right) =E^{-1}\left(
\begin{array}{c}
\phi _{1} \\
\phi _{2}
\end{array}
\right)  \label{psi}
\end{equation}
\end{center}

have the property $\left( \psi _{i},\varphi _{j}\right) =\delta _{ij}.$

By (\ref{psi}) $\psi _{1}\left( s\right) =\frac{2}{2+\pi i}e^{is},\;s\in %
\left[ 0,r\right] ,$ and $\psi _{1}\left( 0\right) =2\frac{2-\pi i}{4+\pi
^{2}}.$

\section{The Bautin-type bifurcation in the central manifold}

We consider $a_{0}=-1$, when the dynamical system admits a two-dimensional
local center manifold, that we denote $W_{loc}^{c}\left( c\right) $.
Everywhere below, the dependence of $\alpha =\left( a,c\right) $ becomes
dependence of $c$ only.

Obviously, we have for every $\phi \in C\left( \left[ -r,0\right] ,\mathbb{R}%
\right) ,$
\begin{eqnarray}
f\left( \left[ S_{c}\left( t\right) \phi \right] \left( s\right) ,\left[
S_{c}\left( t\right) \phi \right] \left( s-r\right) ,c\right) &=&\left[
S_{c}\left( t\right) \phi \right] \left( s\right) \left[ S_{c}\left(
t\right) \phi \right] \left( s\right) +  \label{f} \\
&&+c\left[ S_{c}\left( t\right) \phi \right] \left( s\right) \;\left[
S_{c}\left( t\right) \phi \right] \left( s-r\right) .  \notag
\end{eqnarray}

By writing $w_{c}\left( s,z\left( t\right) ,\overset{\_}{z}\left( t\right)
\right) $ as a series of powers of $z$ and $\overset{\_}{z},$

\begin{equation}
w_{c}\left( s,z,\overset{\_}{z}\right) =\sum\limits_{j+k\geq 2}\frac{1}{j!k!}%
w_{jk}\left( s,c\right) z^{j}\overset{\_}{z}^{k},  \label{w}
\end{equation}

and inserting (\ref{w}) and (\ref{Salfa}) in (\ref{f}), we can obtain the
coefficients $F_{jk}$ of (\ref{F-z}) (here depending on $(s,c)$).

Since

$f\left( S_{c}\left( t\right) \phi \left( 0\right) ,S_{c}\left( t\right)
\phi \left( -r\right) ,c\right) =$

$=\left[ z\left( t\right) \varphi _{1}\left( 0\right) +\overset{\_}{z}\left(
t\right) \overset{\_}{\varphi }_{1}\left( 0\right) +\frac{1}{2}w_{20}\left(
0\right) z^{2}+w_{11}\left( 0\right) z\overset{\_}{z}+\frac{1}{2}%
w_{02}\left( 0\right) \overset{\_}{z}^{2}+...\right] $

$\left[ z\left( t\right) \varphi _{1}\left( 0\right) +\overset{\_}{z}\left(
t\right) \overset{\_}{\varphi }_{1}\left( 0\right) +\frac{1}{2}w_{20}\left(
0\right) z^{2}+w_{11}\left( 0\right) z\overset{\_}{z}+\frac{1}{2}%
w_{02}\left( 0\right) \overset{\_}{z}^{2}+...\right] +$

$+c\left[ z\varphi _{1}\left( 0\right) +\overset{\_}{z}\overset{\_}{\varphi }%
_{1}\left( 0\right) +\frac{1}{2}w_{20}\left( 0\right) z^{2}+w_{11}\left(
0\right) z\overset{\_}{z}+\frac{1}{2}w_{02}\left( 0\right) \overset{\_}{z}%
^{2}+...\right] $

$\left[ z\varphi _{1}\left( -r\right) +\overset{\_}{z}\overset{\_}{\varphi }%
_{1}\left( -r\right) +\frac{1}{2}w_{20}\left( -r\right) z^{2}+w_{11}\left(
-r\right) z\overset{\_}{z}+\frac{1}{2}w_{02}\left( -r\right) \overset{\_}{z}%
^{2}+...\right] ,$

we find, by denoting $F_{jk}\left( 0,c\right) =F_{jk},$

\begin{center}
$\ \ \ \ \ \ \ \ \ \ F_{20}=2\left( 1-ic\right) ,$

$F_{11}=2,$

$\ \ \ \ \ \ \ \ \ \ F_{02}=2\left( 1+ic\right) .$
\end{center}

By the definition of the function $g,$ (\ref{rel f-g}), and of the
coefficients $g_{jk},$

\begin{equation}
g_{jk}=\psi _{1}\left( 0\right) F_{jk}=2\frac{2-\pi i}{4+\pi ^{2}}F_{jk}.
\label{g_F}
\end{equation}

Hence, by the above relations, $g_{20},\;g_{11},\;g_{02}$ are determined.\

Now we look for the second order terms in the series of powers defining $%
w_{c},\;$ (\ref{w}). Differential equations for them are found from the
relation \cite{MNO}, \cite{I}
\begin{eqnarray*}
&&\frac{\partial }{\partial s}\sum\limits_{j+k\geq 2}\frac{1}{j!k!}%
w_{jk}(s,c)z^{j}\overset{\_}{z}^{k}\;=\;\sum\limits_{j+k\geq 2}\frac{1}{j!k!}%
g_{jk}(c)z^{j}\overset{\_}{z}^{k}\varphi _{1}\left( s\right) + \\
&&+\sum\limits_{j+k\geq 2}\frac{1}{j!k!}\overset{\_}{g}_{jk}(c)\overset{\_}{z%
}^{j}z^{k}\overset{\_}{\varphi }_{1}\left( s\right) +\frac{\partial }{%
\partial t}\sum\limits_{j+k\geq 2}\frac{1}{j!k!}w_{jk}(s,c)z^{j}\overset{\_}{%
z}^{k},
\end{eqnarray*}
by equating the terms containing the same powers of $z\left( t\right) $ and $%
\overset{\_}{z}\left( t\right) $. The conditions for the determination of
the integration constants are obtained from \cite{MNO}, \cite{I}
\begin{eqnarray*}
&&\frac{d}{dt}\sum\limits_{j+k\geq 2}\frac{1}{j!k!}w_{jk}(0,c)z^{j}\overset{%
\_}{z}^{k}+\sum\limits_{j+k\geq 2}\frac{1}{j!k!}g_{jk}(c)z^{j}\overset{\_}{z}%
^{k}\varphi _{1}\left( 0\right) +\sum\limits_{j+k\geq 2}\frac{1}{j!k!}%
\overset{\_}{g}_{jk}(c)\overset{\_}{z}^{j}z^{k}\varphi _{2}\left( 0\right)
\;\; \\
&=&-\sum\limits_{j+k\geq 2}\frac{1}{j!k!}w_{jk}(-r,c)z^{j}\overset{\_}{z}%
^{k}+\sum\limits_{j+k\geq 2}\frac{1}{j!k!}F_{jk}z^{j}\overset{\_}{z}^{k}.
\end{eqnarray*}

Thus, we find for $w_{20}\left( s,c\right) ,$

\begin{center}
$w_{20}^{\prime }\;=2iw_{20}(s,c)+g_{20}(c)e^{is}+\overset{\_}{g}%
_{02}(c)e^{-is},$

$2w_{20}\left( 0,c\right) i+g_{20}\left( c\right) +\overset{\_}{g}%
_{02}\left( c\right) =-w_{20}\left( -r,c\right) +2-2ic,\;$
\end{center}

and by solving the equation we have

\begin{center}
$w_{20}(s,c)=w_{20}(0,c)e^{2is}-\frac{1}{i}g_{20}(c)\left(
e^{is}-e^{2is}\right) -\frac{1}{3i}\overset{\_}{g}_{02}(c)\left(
e^{-is}-e^{2is}\right) ,$

$w_{20}\left( 0,c\right) =\frac{2\left( 1+2i\right) }{15\left( 4+\pi
^{2}\right) }\left[ \left( 4-3\pi ^{2}\right) +8c+i\left( 8-4c+3\pi
^{2}c\right) \right] .$
\end{center}

For $w_{11}\left( s,c\right) :$

\begin{center}
$w_{11}^{\prime }=\;g_{11}(c)e^{is}+\overset{\_}{g}_{11}(c)e^{-is},$

$w_{11}\left( -r\right) =w_{11}\left( 0\right) +ig_{11}(\alpha )\left(
i+1\right) +i\overset{\_}{g}_{11}(\alpha )\left( i-1\right) ,\;\;\;$
\end{center}

from where

\begin{center}
$w_{11}\left( s\right) =w_{11}\left( 0\right) -ig_{11}(c)\left(
e^{is}-1\right) +i\overset{\_}{g}_{11}(c)\left( e^{-is}-1\right) ,$

$w_{11}\left( 0\right) =2\frac{\pi ^{2}-4}{4+\pi ^{2}}+g_{11}(c)\left(
1-i\right) +\overset{\_}{g}_{11}(c)\left( 1+i\right) .$ \ \ \ \ \ \ \ \ \ \
\ \ \
\end{center}

The relation $w_{02}=\overset{\_}{w}_{20}$holds true.

\subsection{The first Lyapunov coefficient}

We are now able to calculate $F_{jk}$ (and thus $g_{jk}$) with $j+k=3$.

We find

\begin{center}
$F_{30}=3c\left( w_{20}\left( -r\right) -w_{20}\left( 0\right) i\right)
+6w_{20}\left( 0\right) ,\;\;\;\;\;\;\;\;\;F_{03}=\overset{\_}{F}%
_{30},\;\;\;\;\;\;\;\;\;\;\;\;\;$

$F_{21}=2c\left( w_{11}\left( -r\right) -iw_{11}\left( 0\right) \right)
+4w_{11}\left( 0\right) +c\left( w_{20}\left( -r\right) +w_{20}\left(
0\right) i\right) +2w_{20}\left( 0\right) ,$

$F_{12}=\overset{\_}{F}_{21},\;$ \ \ \ \ \ \ \ \ \ \ \ \ \ \ \ \ \ \ \ \ \ \
\ \ \ \ \ \ \ \ \ \ \ \ \ \ \ \ \ \ \ \ \ \ \ \ \ \ \ \ \ \ \ \ \ \ \ \ \ \
\end{center}

while $g_{jk}$ are given by (\ref{g_F}).

These allow us to calculate the first Lyapunov coefficient, (\ref{l1}):
\begin{equation}
l_{1}\left( c\right) =\frac{1}{5(4+\pi ^{2})}\left[ \left( 8-12\pi \right)
c^{2}+\left( 72-28\pi \right) c+144-16\pi \right] .  \label{L1}
\end{equation}

We impose the condition

\begin{center}
\bigskip $l_{1}\left( c\right) =0\;\;\Leftrightarrow \;\;\;\left( 8-12\pi
\right) c^{2}+\left( 72-28\pi \right) c+144-16\pi =0.$
\end{center}

The two solutions of this equations are
\begin{eqnarray*}
c_{1} &=&\frac{18-7\pi +\sqrt{36+212\pi +\pi ^{2}}}{2(3\pi -2)}\thickapprox
1.52799, \\
c_{2} &=&\frac{18-7\pi -\sqrt{36+212\pi +\pi ^{2}}}{2(3\pi -2)}\approx
-2.06554.
\end{eqnarray*}

We thus found two values of the parameter $c$ for which degenerate Hopf
bifurcation takes place.

\subsection{The second Lyapunov coefficient}

The second Lyapunov coefficient at the values $c_{1}$, $c_{2}$ of the
parameter $c$ has the form \cite{K}

\noindent $12l_{2}(c_{i})=\frac{1}{\omega _{0}}\mathrm{Re}g_{32}+\;\;\;\;\;\;%
\;\;\;\;\;\;\;\;\;\;\;\;\;\;\;\;\;\;\;\;\;\;\;\;\;\;\;\;\;\;\;\;\;\;\;\;\;\;%
\;\;\;\;\;\;\;\;\;\;\;\;\;\;\;\;\;\;\;\;\;\;\;\;\;\;\;$

$\;\;\;\;\ \ \ +\frac{1}{\omega _{0}^{2}}\mathrm{Im}\left[ g_{20}\overset{\_}{g%
}_{31}-g_{11}\left( 4g_{31}+3\overset{\_}{g}_{22}\right) -\frac{1}{3}%
g_{02}\left( g_{40}+\overset{\_}{g}_{13}\right) -g_{30}g_{12}\right]
+\;\;\;\;\;\;$

$\;\;\;\ \ \ \ +\frac{1}{\omega _{0}^{3}}\{\mathrm{Re}\left[ g_{20}\left(
\overset{\_}{g}_{11}\left( 3g_{12}-\overset{\_}{g}_{30}\right) +g_{02}\left(
\overset{\_}{g}_{12}-\frac{1}{3}g_{30}\right) +\frac{1}{3}\overset{\_}{g}%
_{02}g_{03}\right) \right. \;\;\;\;\;\;\;\;\;\;$

$\;\;\;\;\;\ \left. +g_{11}\left( \overset{\_}{g}_{02}\left( \frac{5}{3}%
\overset{\_}{g}_{30}+3g_{12}\right) +\frac{1}{3}g_{02}\overset{\_}{g}%
_{03}-4g_{11}g_{30}\right) \right] +3\mathrm{Im}\left( g_{20}g_{11}\right)
\mathrm{Im}g_{21}\}+\;$

$+\frac{1}{\omega _{0}^{4}}\{\mathrm{Im}[g_{11}\overset{\_}{g}_{02}\left(
\overset{\_}{g}_{20}^{2}-3\overset{\_}{g}_{20}g_{11}-4g_{11}^{2}\right) ]+%
\mathrm{Im}\left( g_{20}g_{11}\right) [3\mathrm{Re}\left( g_{20}g_{11}\right)
-2\left| g_{02}\right| ^{2}]\},$

where $g_{ij}$ are evaluated at $c_{i},$ $i=1$ or $i=2.$

We will calculate first $l_{2}(c_{1}),$ and for this, in the sequel,\ all
the $g_{jk}$ will be evaluated at $c_{1}.$

Hence, we have to calculate $g_{jk}\left( c_{1}\right) $ for $j+k=4$ and $%
g_{32}.$ In order to do this, we calculate $w_{jk},$ with $j+k=3.$ For $%
w_{30},$ the equation is

\begin{center}
$w_{30}^{\prime }\left( s\right) =3iw_{30}\left( s\right) +g_{30}e^{is}+%
\overset{\_}{g}_{03}e^{-is}+$

$+3w_{20}(s)g_{20}+3w_{11}(s)\overset{\_}{g}_{02},$
\end{center}

\bigskip \noindent with the condition

\begin{center}
$w_{30}\left( -r\right) =-3w_{30}\left( 0\right) i-3w_{20}\left( 0\right)
g_{20}-3w_{11}\left( 0\right) \overset{\_}{g}_{02}-g_{30}-\overset{\_}{g}%
_{03}+$

$+3c_{1}\left( w_{20}\left( -r\right) -w_{20}\left( 0\right) i\right)
+6w_{20}\left( 0\right) .$
\end{center}

The values obtained after solving the above equtions, are

\begin{center}
$w_{30}\left( 0\right) =.327626-5.115802i,\;\;\;\;\;w_{03}\left( 0\right) =%
\overset{\_}{w}_{30}\left( 0\right) ,$

$w_{30}\left( -r\right) =-14.190120-5.277852i,\;w_{03}\left( -r\right) =%
\overset{\_}{w}_{30}\left( -r\right) .$
\end{center}

For $w_{21}$, we have

\begin{center}
\bigskip $w_{21}^{\prime }\left( s\right) \;=iw_{21}+g_{21}e^{is}+\overset{\_%
}{g}_{12}e^{-is}+2w_{20}(s)g_{11}+$

$+w_{11}\left( s\right) \left( g_{20}+2\overset{\_}{g}_{11}\right) +w_{02}(s)%
\overset{\_}{g}_{20}$
\end{center}

and

\begin{eqnarray}
w_{21}\left( -r\right) +iw_{21}\left( 0\right) &=&-2w_{20}\left( 0\right)
g_{11}-2w_{11}\left( 0\right) \overset{\_}{g}_{11}-w_{11}\left( 0\right)
g_{20}  \label{cond1} \\
&&-w_{02}\left( 0\right) \overset{\_}{g}_{02}-g_{21}-\overset{\_}{g}%
_{12}+F_{21}.  \notag
\end{eqnarray}

In this case, after solving the differential equation above, the second
equation for obtaining $w_{21}\left( -r\right) $ and $w_{21}\left( 0\right) $
reads
\begin{eqnarray}
w_{21}\left( -r\right) +iw_{21}\left( 0\right) &=&F_{21}-F_{21}\frac{8}{%
4+\pi ^{2}}+2g_{11}e^{-ir}\int_{0}^{-r}w_{20}(\tau )e^{-i\tau }d\tau
\label{cond2} \\
&&+e^{-ir}\left[ \left( g_{20}+2\overset{\_}{g}_{11}\right)
\int_{0}^{-r}w_{11}\left( \tau \right) e^{-i\tau }d\tau +\overset{\_}{g}%
_{02}\int_{0}^{-r}w_{02}(\tau )e^{-i\tau }d\tau \right] .  \notag
\end{eqnarray}
The system (\ref{cond1}), (\ref{cond2}) is not determined. It is to be seen
whether it is compatible or not.

We have to compare the right hand sides of (\ref{cond1}) and (\ref{cond2}).
We first notice that

\begin{center}
$g_{21}+\overset{\_}{g}_{12}=\frac{8}{4+\pi ^{2}}F_{21}.$
\end{center}

Then, we find by direct calculations that
\begin{eqnarray*}
I_{1} &=&\int_{0}^{-r}w_{20}(\tau )e^{-i\tau }d\tau
=w_{20}(-r)+iw_{20}(0)-irg_{20}+\overset{\_}{g}_{02}, \\
I_{2} &=&\int_{0}^{-r}w_{11}(\tau )e^{-i\tau }d\tau
=-w_{11}(-r)-iw_{11}(0)+irg_{11}-\overset{\_}{g}_{11}, \\
I_{3} &=&\int_{0}^{-r}w_{02}(\tau )e^{-i\tau }d\tau =\frac{1}{3}\left[
-w_{02}(-r)+iw_{02}(0)+irg_{02}-\overset{\_}{g}_{20}\right] ,
\end{eqnarray*}
and, by using these equalities,

\begin{center}
$\;\;\;\;\;\;\;\;\;\;\;\;\;\;\;\;\;\;\;\;\;\;-2ig_{11}I_{1}+2w_{20}\left(
0\right) g_{11}=0,\;\;\;\;$

$-i\left[ g_{20}+2\overset{\_}{g}_{11}\right] I_{2}+\left( 2\overset{\_}{g}%
_{11}+g_{20}\right) w_{11}\left( 0\right) =0,$

$\;\;\;\;\;\;\;\;\;\;\;\;\;\;\;\;-i\overset{\_}{g}_{02}I_{3}+w_{02}\left(
0\right) \overset{\_}{g}_{02}=0.$
\end{center}

This leads to the conclusion that the right hand sides of (\ref{cond1}) and (%
\ref{cond2}) are equal, and thus, the system in $w_{21}\left( -r\right) $, $%
w_{21}\left( 0\right) $ is compatible. We will then take $w_{21}\left(
0\right) =0$ and it follows that

\begin{center}
$w_{21}\left( -r\right) =-\;2w_{20}\left( 0\right) g_{11}-2w_{11}\left(
0\right) \overset{\_}{g}_{11}-w_{11}\left( 0\right) g_{20}$

$\;\;-w_{02}\left( 0\right) \overset{\_}{g}_{02}-g_{21}-\overset{\_}{g}%
_{12}+F_{21}.$
\end{center}

Then $w_{12}\left( 0\right) =0$ and $w_{12}\left( -r\right) =\overset{\_}{w}%
_{21}\left( -r\right) .$

We are now able to compute $g_{jk},$ with $j+k=4.$

Firstly we compute $F_{jk},$ with $j+k=4:$

$\frac{1}{24}F_{40}=\frac{1}{3}w_{30}\left( 0\right) +\frac{1}{4}%
w_{20}\left( 0\right) ^{2}+c_{1}\left( \frac{1}{6}w_{30}\left( -r\right) -%
\frac{1}{6}iw_{30}\left( 0\right) +\frac{1}{2}w_{20}\left( 0\right) \frac{1}{%
2}w_{20}\left( -r\right) \right) ,$

$\left. {}\right. $

$\frac{1}{6}F_{31}\;=c_{1}\left( \frac{1}{2}w_{21}\left( -r\right) +\frac{1}{%
6}w_{30}\left( -r\right) \right. +\frac{1}{6}iw_{30}\left( 0\right) +\left. +%
\frac{1}{2}w_{20}\left( 0\right) w_{11}(-r)\right) +$

$\ \ \ \ \ \ \ \ \ \ \ \;+w_{11}(0)\frac{1}{2}w_{20}\left( 0\right) +\frac{1%
}{3}w_{30}\left( 0\right) +w_{20}\left( 0\right) w_{11}(0),$

$\left. {}\right. $

$\frac{1}{4}F_{22}\;=c_{1}\left[ \frac{1}{2}w_{12}\left( -r\right) +\frac{1}{%
2}w_{21}\left( -r\right) \right. +\frac{1}{4}w_{20}\left( 0\right)
w_{02}\left( -r\right) +w_{11}(0)w_{11}(-r)+$

$\;\;\;\;\;\;\;\;\left. +\frac{1}{4}w_{02}\left( 0\right) w_{20}\left(
-r\right) \right] +\left[ \frac{1}{2}w_{20}\left( 0\right) w_{02}\left(
0\right) +w_{11}(0)w_{11}(0)\right] .$

\bigskip

We obtain, by using (\ref{g_F}):

\bigskip

$g_{40}=-70.452908+32.020324i,$ $\;\;g_{04}=.804019+77.383894i,$

$g_{31}=-13.491939-6.450063i,\;\;\;\;g_{13}=11.553771+9.494531i,$

$g_{22}=4.485812-7.046298i.$

$\left. {}\right. $

The only $g_{jk}$ still to be computed in order to be able to evaluate $%
l_{2}(c_{1})$ is $g_{32}.\;$Since

\begin{center}
$F_{32}=6w_{22}\left( 0\right) +4w_{31}\left( 0\right) +6w_{20}\left(
0\right) w_{12}\left( 0\right) +12w_{11}\left( 0\right) w_{21}\left(
0\right) +2w_{02}\left( 0\right) w_{30}\left( 0\right) +$

$+c\left[ 3w_{22}\left( -r\right) +2w_{31}\left( -r\right) +\right.
+3w_{20}\left( 0\right) w_{12}\left( -r\right) +6w_{11}\left( 0\right)
w_{21}\left( -r\right) +$

$\ \ \ \ \ \ \ \ \ \ \ +w_{02}\left( 0\right) w_{30}\left( -r\right)
+w_{30}\left( 0\right) w_{02}\left( -r\right) +6w_{21}\left( 0\right)
w_{11}\left( -r\right) +\;\;\;\;\;\;\;\;\;\;\;\;\;\;\;\;\;\;\;$

$\;\;\;\;\;\;\;\;\;\;\;\;\;\;\;+3w_{12}\left( 0\right) w_{20}\left(
-r\right) +2w_{31}\left( 0\right) i+3w_{22}\left( 0\right) \left( -i\right)
],$ \ \ \ \ \ \ \ \ \ \ \ \ \ \ \ \ \ \ \ \ \ \ \ \ \ \ \ \ \ \ \ \ \ \ \ \
\ \ \ \ \ \ \ \ \ \ \ \ \ \ \ \

\ \ \ \ \ \ \ \ \ \ \ \ \ \ \ \ \ \ \ \ \ \ \ \ \ \ \ \ \ \ \ \ \ \ \ \ \ \
\ \ \ \ \ \ \ \ \ \ \ \ \ \ \ \ \ \
\end{center}

\noindent we have to compute $w_{22}(0),\;w_{22}(-r),\;w_{31}(0),%
\;w_{31}(-r).$

The equations for $w_{22}(s)$ are

\begin{center}
$\;w_{22}^{\prime }=g_{22}e^{is}+\overset{\_}{g}%
_{22}e^{-is}+2w_{20}(s)g_{12}+2w_{02}(s)\overset{\_}{g}_{12}\;\;\;\;\;\;\;\;%
\;\;\;\;\;\;\;\;\;\;\;\;\;\;\;$

$+2w_{11}(s)\left( g_{21}+\overset{\_}{g}_{21}\right)
+w_{30}(s)g_{02}+w_{03}(s)\overset{\_}{g}_{02}+\;\;\;\;\;\;\;\;\;\;\;$

$\ \ \ \ \ \ \ \ \ \ \ \ +w_{21}(s)\left( 4g_{11}+\overset{\_}{g}%
_{20}\right) +w_{12}(s)\left( g_{20}+4\overset{\_}{g}_{11}\right)
,\;\;\;\;\;\;\;\;\;\;\;\;\;\;\;\;\;\;\;\;\;\;\;\;\;$
\end{center}

and$\;\;\;\;\;\;\;$

\begin{center}
$w_{22}\left( -r\right) =-(2w_{20}\left( 0\right) g_{12}+2w_{11}\left(
0\right) \overset{\_}{g}_{21}+2w_{11}\left( 0\right) g_{21}++2w_{02}\left(
0\right) \overset{\_}{g}_{12}+\;\;\;\;\;\;\;$

$+w_{30}\left( 0\right) g_{02}+4w_{21}\left( 0\right) g_{11}++w_{21}\left(
0\right) \overset{\_}{g}_{20}+w_{12}\left( 0\right) g_{20}+$

$+4w_{12}\left( 0\right) \overset{\_}{g}_{11}+w_{03}\left( 0\right) \overset{%
\_}{g}_{02}+g_{22}+\overset{\_}{g}_{22})+\;\;\;\;\;\;\;\;\;\;\;$

$+4w_{12}\left( 0\right) +4w_{21}\left( 0\right) +2w_{20}\left( 0\right)
w_{02}\left( 0\right) +4w_{11}(0)w_{11}(0)+$

$+c[2w_{12}\left( -r\right) +2w_{21}\left( -r\right) +2\left( -i\right)
w_{12}\left( 0\right) +2iw_{21}\left( 0\right) +\;\;\;\;\;\;$

$+w_{20}\left( 0\right) w_{02}\left( -r\right)
+4w_{11}(0)w_{11}(-r)+w_{02}\left( 0\right) w_{20}\left( -r\right)
]\;\;\;\;\;\;$
\end{center}

while those for $w_{31}(s)$ are

\begin{center}
$w_{31}^{\prime }=2w_{31}i+g_{31}e^{is}+\overset{\_}{g}%
_{13}e^{-is}+3w_{20}(s)g_{21}+w_{11}(s)g_{30}+$

$+3w_{11}(s)\overset{\_}{g}_{12}+w_{02}(s)\overset{\_}{g}_{03}+3w_{30}\left(
s\right) g_{11}+$

$\;+3w_{21}\left( s\right) g_{20}+3w_{21}\left( s\right) \overset{\_}{g}%
_{11}+3w_{12}\left( s\right) \overset{\_}{g}_{02},$
\end{center}

and

\begin{center}
$2w_{31}\left( 0\right) i+w_{31}\left( -r\right) =-(3w_{20}\left( 0\right)
g_{21}(\alpha )+3w_{11}\left( 0\right) \overset{\_}{g}_{12}(\alpha
)+w_{11}\left( 0\right) g_{30}(\alpha )+$

$\;\;\;\;\;\;\;\;\;\;\;\;+w_{02}\left( 0\right) \overset{\_}{g}%
_{03}+3w_{30}\left( 0\right) g_{11}+3w_{21}\left( 0\right) g_{20}+$

$\;\;\;\;\;\;\;\;\;\;\;\;\;\;\;\;\;+3w_{21}\left( 0\right) \overset{\_}{g}%
_{11}+3w_{12}\left( 0\right) \overset{\_}{g}_{02}+g_{31}+3\overset{\_}{g}%
_{13})+$

$\;\;\;\;\;\;\;\;\;\;\;\;\;\;+c[3w_{21}\left( -r\right) +w_{30}\left(
-r\right) -3iw_{21}\left( 0\right) +iw_{30}\left( 0\right) +$

$\;\;\;\;\;\;+3w_{20}\left( 0\right) w_{11}(-r)+3w_{11}(0)w_{20}\left(
-r\right) ]+$

$\;\;\;\;\;\;\;\ +\left[ 6w_{21}\left( 0\right) +2w_{30}\left( 0\right)
+6w_{20}\left( 0\right) w_{11}(0)\right] .$
\end{center}

After computations:

\begin{center}
$w_{22}(-r)=4.864870928,$

$w_{22}(0)=-43.85187247,$

$w_{31}\left( -r\right) =-6.41714235-18.89415271i,$

$w_{31}\left( 0\right) =17.94690049+2.001612024i.$
\end{center}

For $g_{32\text{ }}$we found the value

\begin{center}
$g_{32}=28.68605342+128.6141166i.$
\end{center}

Now we can compute $l_{2}(c_{1}).$ We find:
\begin{equation*}
l_{2}(c_{1})=13.08553919.
\end{equation*}

Hence $l_{2}(c_{1})>0.$ We are now able to assert and prove

\textbf{Proposition 2. }\textit{Hypothesis H2 is satisfied by equation }(\ref
{mainec}).

\textbf{Proof }The only part of Hypothesis \textit{H2} that still has to be
checked, is that the map $\left( a,c\right) \rightarrow \left( \nu _{1},\nu
_{2}\right) \,\ $is regular at $(-1,c_{1})$, where $\nu _{1}=\frac{\mu
\left( a\right) }{\omega \left( a\right) },\;\nu _{2}=l_{1}\left( a,c\right)
.$ We have
\begin{equation*}
\frac{\partial \left( \nu _{1},\nu _{2}\right) }{\partial \left( a,c\right) }%
=\left(
\begin{array}{cc}
\left( \frac{\mu \left( a\right) }{\omega \left( a\right) }\right) ^{\prime }
& 0 \\
\frac{\partial }{\partial a}l_{1}\left( a,c\right) & \frac{\partial }{%
\partial c}l_{1}\left( a,c\right)
\end{array}
\right) .
\end{equation*}
We have $\left. \left( \frac{\mu \left( a\right) }{\omega \left( a\right) }%
\right) ^{\prime }\right| _{a=-1}=\left. \frac{\mu ^{\prime }\left( a\right)
\omega \left( a\right) -\mu \left( a\right) \omega ^{\prime }\left( a\right)
}{\omega ^{2}\left( a\right) }\right| _{a=-1}=\frac{\mu ^{\prime }\left(
-1\right) }{\omega \left( -1\right) }=\mu ^{\prime }\left( -1\right) .$ By
taking the derivative of the equation for $\mu $, (\ref{ec_cos}), with
respect to $a$, and by evaluating the result in $a=-1,$ we find $\mu
^{\prime }\left( -1\right) \neq 0.$

The form of $l_{1},$ and the fact that the equation $l_{1}(-1,c)=0$ has two
distinct solutions show that $\frac{\partial }{\partial c}l_{1}\left(
-1,c_{1}\right) \neq 0.$

It follows that $\left. \frac{\partial \left( \nu _{1},\nu _{2}\right) }{%
\partial \left( a,c\right) }\right| _{(-1,c_{1})}\neq 0$ hence the
conclusion.$\square $

Since all the hypotheses of the \textbf{Theorem }presented in Introduction
are satisfied, we may formulate the following result.

\bigskip \textbf{Proposition 3. }\textit{The equation (\ref{mainec})
presents a Bautin-type bifurcation at}

\noindent $\left( a_{0},c_{0}\right) =(-1,c_{1}).$

\bigskip

Following the same path as for $l_{1}\left( c_{1}\right) $, we find $%
l_{2}(c_{2})<0.$ In this situation, as is shown in \cite{K}, the two limit
cycles (one interior to the other) that should appear by the Bautin
bifurcation for eq. (\ref{ecz}), exist for some zone of the quadrant $\nu
_{1}<0,\;\nu _{2}>0.$ But $\nu _{1}<0\;\Leftrightarrow \;\mu (a,c)<0,$ and
we have no theorem to assert the existence of a bi-dimensional invariant
(stable) manifold that is tangent to $\mathbb{M}_{\left\{ \lambda
_{1,2}\left( \alpha \right) \right\} }.$ That is why the restriction $%
l_{2}(\alpha )>0$ is among the hypotheses of our \textbf{Theorem}.

\end{document}